\documentclass[reqno]{amsart}
\usepackage{amsfonts,amsthm,amsbsy,amsmath,amssymb}
\usepackage{latexsym,enumerate}
\usepackage{wasysym}

\newtheorem{corollary}{Corollary}
\newtheorem{lemma}{Lemma}
\newtheorem{prop}{Proposition}
\newtheorem{theorem}{Theorem}
\newtheorem{remark}{Remark} 

\def\C{{\mathbb C}}
\def\N{{\mathbb N}}

\def\R{{\mathbb R}}

\def\Cinf{{\mathbb C_\infty}}

\def\a{\alpha}

\def\e{\varepsilon}
\def\g{\gamma}
\def\s{\sigma}

\begin{document}
\title[Connectedness of the complement of a union of two compact sets
in $\R^2$]{The Union of two compact sets in $\R ^2$ with connected complement has a connected complement}

\author[Y. Dossidis]{Yeorgios A. Dossidis}

\email{gntosi01@ucy.ac.cy}

\begin{abstract}
In the papers from Chui and Parnes (1971) \cite{CP} and Luh (1972) \cite{Luh}, as well on the paper from V.Nestoridis (1996) \cite{N} on the Universal Taylor series, it is used, without proof, that the union of two compact sets in $\R ^2$ with connected complement has a connected complement. In this work we present a rigorous proof of this fact, by studying some important topological properties of the plane.
\end{abstract}

\maketitle
 
\section*{Preliminaries}

The following known results will be used in proving the main theorem.

\begin{prop}\label{prop-one}
Let $A\subset\C$ open. Then $A$ is path-connected if and only if $A$ is connected. 

\proof \noindent Every path-connected subset of $\C$ is connected, so it suffices to show 
the other direction.
Let $A\subset\C$ be open and connected. That implies that if the exist
$U,V\subset A$ open, such as $U\cup V= A$ and $U\cap V = \varnothing$, 
then either $U=\varnothing$ or $V=\varnothing$.
Assume $A\ne\varnothing$ and $x_0 \in A$, and define 
\begin{align*} 
U:=\{x\in A \text{ } | \text{ } \exists C:[0,1]\to A \text{, continuous, such as } C(0)=x_0 
\text{and\,\,}C(1)=x\}. 
\end{align*}
Then $U$ is open. Indeed, if $x\in U$ then there exists an 
$\e > 0$, such that $B(x,\e)\subset A$ because $A$ is open.
Now, for all $x_1\in B(x,\e)$, there exists $C_1:[0,1] \to A$, 
such as $C_1(0)= x$, $C_1(1)=x_1$ continuous (the line segment $xx_1$), and so 
$$
C_2:=C\oplus C_1 = \left\{
     \begin{array}{lcr}
       C(2s) & \text{if} & s \in [0,\frac{1}{2}], \\
       C_1(2s-1) & \text{if} & s \in \big[\frac{1}{2},1\big],
     \end{array}
   \right.
$$
is continuous and  $C_2(0)=x_0$, $C_2(1)=x_1$. So $x_1 \in U$ which implies $B(x,\e)\subset U$ 
and so $U$ is open. However, $V:=A\setminus U$ is also open: If $x \in V$ then there exists an
$\e_2>0 : B(x,\e_2)\subset A$. If there exists an
$x_2\in B(x,\e_2):x_2\in U$ then $x\in U$, which leads to a contradiction.
So $B(x,\e_2)\subset V$ which implies $V$ open.
As $U,V$ are open and $U\cup V=A$, $U\cap V= \varnothing$ 
it follows that $V=\varnothing$ and $A=U$, which means $A$ is path-wise connected. \hfill $\Box$

\end{prop}

\begin{remark}\label{rem_one}
The proof above can be slightly modified in order to be valid for topological
space which is locally path-connected.
In particular, Proposition \ref{prop-one} holds for any open $A\subset\Cinf=\C\cup\{\infty\}$, where 
$\Cinf$ the Riemann sphere.
\end{remark}

\begin{prop}\label{prop_two}
Let $A\subset\C$ (or $\Cinf$) bounded. Then $A^c$ has exactly one unbounded component.  

\proof Because $A$ is bounded, then there exists 
$M >0 : A\subset B(0,M) \Rightarrow A^c \supset \big(B(0,M)\big)^c$ and so $A^c$ 
is not bounded. Furthermore, $\big(B(0,M)\big)^c$ is connected and $A^c \setminus \big(B(0,M)\big)^c$ 
is bounded or empty. Therefore $A^c$ has exactly one unbounded component and every other 
component (if any exist) are contained in $B(0,M)$ \hfill $\Box$ \\

\end {prop}

\begin{prop}\label{prop_three}
\noindent Let $A\subset\C$ bounded. Then $\C \setminus A$ is connected if and only if 
$\Cinf \setminus A $ is connected. 
\proof Let $A$ bounded. Then $\C \setminus A$ has exactly one unbounded component, $U_{A^c}$.\\If $\C \setminus A$ is connected, then $\C \setminus A = U_{A^c}$ and therefore every neighborhood of infinity is a subset of $U_{A^c}\cup \{ \infty\}$. We conclude  that $\Cinf \setminus A = U_{A^c}\cup \{ \infty\} $ is connected.\\
Lets now assume $\C \setminus A$ disconnected. Then, because $A$ is bounded, there exists $B \neq \varnothing$, a bounded component of $\C \setminus A$  ( $B\subset (\C \setminus A)\setminus U_{A^c}$).\\ Then $B\cap (U_{A^c}\cup \{\infty\}) = \varnothing$ and therefore $\Cinf \setminus A$ disconnected. \hfill $\Box$ 
\end{prop}

\section*{Simply Connected Subsets of the Plane}

\noindent The following theorem is a characterization of simply connected subsets of the plane. 

\begin{theorem}\label{theorem_one}
An open and connected $A\subset\C$ is simply connected if and only if $\Cinf \setminus A$ is connected.  
\\
\noindent In order to prove this, we will need the following \footnote{I am indebted to Y.Smyrlis for suggesting this line of proof.}:

\begin{lemma}\label{lemma1}
Let $X$ be a connected topological manifold, and $x,y \in X$ such as $x\neq y$. Then there exists an injective curve $\s : [0,1] \to X$ with $\s(0)=x$ and $\s(1)=y$
\proof (of lemma \ref{lemma1})\\
Let $x_0 \in X$ and define $W=\{x\in X: x\neq x_0$ and $x_0$ is connected to $x$ by an injective curve$\}\cup\{x_0\}$.\\
Now $W$ is both open and closed. Indeed, let $x_1\in W$, and $V$ an open neighbourhood of $x_1$ which is homeomorphic to the unit ball $B(0,1)$, with $p: V\to B(0,1)$ a homeomorphism with $p(x_1)=0$\\
Let $x_2\in V$. Without any loss of generality, we can assume that $x_2$ does not belong to the injective curve $\s_1$ which connects $x_0$ and $x_1$. If $r:=\|p(x_2)\|<1$,  let \begin{align*}d=\mathrm{dist}\big(p(x_2),\overline{B}(0,r)\cap p(\gamma[I])\big)>0.
\end{align*}
Let $y_3\in \overline{B}(0,r)$, such that $\|y_3-p(x_2)\|=d$ (such $y_3$ exists since $\overline{B}(0,r)\cap p(\gamma[I])$ is compact), $J$ be the segment connecting $y_3$ and $p(x_2)$ and $x_3=p ^{-1}(y_3)$. Then the union of $p^{-1}(J)$ and the part of $\s_1$ which connects $x_0$ with $x_3$ is an injective curve connecting $x_0$ and $x_2$, and therefore $W$ is open.\\
Using a similar approach to the one we just used, it is readily obtained that a limit point of $W$ also belong to $W$.\\
Since $X$ is assumed connected, and $W\neq \varnothing$, we have that $W = X$ and, therefore such an injective curve exists for any $x,y \in X$ with $x\neq y$ \hfill $\Box$
\end{lemma}

\begin{prop}\label{prop A1}
Let $A \subset \C$ be open, connected and not simply connected. Then there exists in $A$ a simple closed curve, which is not $A$-null-homotopic.
\proof
Let $X$ be the universal cover of $A$ and $\pi : X \to A$ be the universal covering map.\\
Since $A$ is not null-homotopic, there exists (by definition) a closed curve $\g: [0,1] \to A$ which is not $A$-null-homotopic.\\ 
Let $\hat{\g} : [0,1] \to X$ be a lift of $\g$. Since $\g$ is not $A$-null-homotopic, we have $\hat{\g}(1) \neq \hat{\g}(0)$ \cite{Forster}. Since $X$ is a connected topological manifold,
%reference (px forster Lectures in Riemann Surfaces)
by the previous lemma, there exists an injective curve $\s :[0,1] \to X$ such that $\s(0)= \hat{\g}(0)$ and $\s(1)= \hat{\g}(1)$. Next define \begin{align*}
s = \inf \bigl\{ t \in (0,1] : \bigl(\exists u \in [0,t)\bigr) : \bigl(\pi(\s(t)) = \pi(\s(u))\bigr)\bigr\}.
\end{align*}
We claim that there is an $u \in [0,s)$ with $\pi(\s(s)) = \pi(\s(u))$: Choose a neighbourhood $V$ of $\s(s)$ on which $\pi$ is injective. There is an $\e > 0$ such that $\s(t) \in V$ for all $t \in [s-\e,s+\e]$. We have sequences $t_n$ and $u_n$ with $\pi(\s(t_n)) = \pi(\s(u_n))$ where $t_n \geqslant t_{n+1} \geqslant s$ and $u_n < t_n$. Without loss of generality, we can assume $t_n < s + \e$ for all $n$. Since $\pi$ is injective on $V$, we have $u_n < s - \e$ for all $n$. After taking a subsequence, we can assume that $u_n$ converges to some $u \leqslant s - \e$. By continuity,

\begin{align*}\pi(\s(s)) = \lim_{n\to\infty} \pi(\s(t_n)) = \lim_{n\to \infty} \pi(\s(u_n)) = \pi(\s(u)).\end{align*}

Now define $\tilde{\s}(t) = \s((1-t)u + ts)$ and $\tilde{\g} = \pi \circ \tilde{\s}$. By construction, $\tilde{\g}$ is a simple closed curve in $U$, and since $\s(u) \neq \s(s)$, it is not null-homotopic. \hfill $\Box$

\end{prop}

\newpage

We are now ready to prove theorem \ref{theorem_one}

\proof of Theorem \ref{theorem_one}\\
 Let A be not simply connected. Then there exists $\g:[0,1] \to A$ simple, continuous and closed that is not $A$-null-homotopic. \\ Let $\g^*:=\g([0,1])$. From the Jordan-Schoenflies curve theorem, we have that $\Cinf\setminus\g^*$ consists of exactly two components, the interior $int(\g)$, which is the bounded one, and the exterior, $ext(\g)$, both of which are homeomorphic to the unit disk.\\ Therefore, if $int(\g) \subset A$, it would follow that $\g$ is $A$-null-homotopic. Therefore $\exists z_1 \in \C : z_1 \in int(\g)$ and $z_1 \not \in A$. Then, on the Riemann sphere, $z_1$ and $\infty$ belong on different components of $\Cinf \setminus \g^*$ and hence $\Cinf \setminus A$ is disconnected.\\
\\
Conversely, let $\Cinf \setminus A$ be disconnected. Let $K$ be a component of $\Cinf \setminus A$ such that $\infty \not \in K$. Then $K$ is compact. Moreover, let $\Omega := \C \setminus U_{A^c}$. Since $ U_{A^c}$ is closed, $\Omega$ is open and since $K \subset \Omega$, there exists $\e > 0$ such that the $dist(K,U_{A^c}) = 2\e$. \\If $z \in K$ then there exists a continuous and closed curve, $\a : [0,1] \to A$, such as and $z$ belongs in a bounded component of $\C \setminus \a^*$. A detail construction of such a curve can be found in \cite{Rudin}. We begin by constructing a grid of horizontal and vertical lines such as the distance between any two adjacent horizontal lines and any two adjacent vertical lines is $\e$. From the squares (closed 2-cells) formed by the grid, we can choose those that intersect with $K$. Since $K$ is compact, only finitely many are required. By parameterizing the edges of these squares carefully, we arrive at a representative of an element of the fundamental group of $A$ which is different from the identity.  \hfill $\Box$
\end{theorem}
\hfill \\
\begin{corollary}\label{cor_one}
Let $A\subset \C$ open. Then every component of $A$ is simply connected if and only $\Cinf \setminus A$ connected.
\proof First lets assume that $\Cinf \setminus A$ is connected. If $A_0 \subset A$ be one connected component of $A$, it suffices to show that $\Cinf \setminus A_0$ connected. We can assume that $A_0 \neq A$ (this is the case of theorem \ref{theorem_one}) and so let $z \in A \setminus A_0$. If $z_0 \in A_0$, then if $\g:[0,1]\to \C $ is continuous such as $\g(0)=z_0$ and $\g(1)=z$ we have $\g^* \not\subset A $ (where $\g^*:= \g([0,1])$ the image of $\g$). Therefore $\exists s_1 \in (0,1)$ and $z_1 \not\in A$ such as $\g(s_1)=z_1$ and consequently $z$ belongs in the only component of $\Cinf \setminus A_0$ which is therefore connected. Now, from theorem \ref{theorem_one}, $A_0$ is simply connected.\\
Let now assume that $\Cinf \setminus A$ disconnected. So, and because the connected components of $A$ are pairwise disjoint, $\exists z \in \C$ which is contained in the bounded component of the compliment of at least one connected component of A. Let $A_1$ be one such component of $A$. Then $\Cinf \setminus A_1$ disconnected and therefore, from theorem \ref{theorem_one}, $A_1$ is not simply connected. \hfill $\Box$ 
\end{corollary}

\newpage
\section*{Main Theorem}

\noindent Having the results above established, we can now prove the main theorem.

\begin{theorem}[Characterazation]\label{main}
Let $K \subset \C$ compact. Then $\C \setminus K$ connected if and only if for all open $A\supset K$ there exists $V$ open such as $K\subset V \subset A$ and every component of $V$ is simply connected.
\proof Let $K^c:=\C \setminus K$ connected. Since $K$ is compact, $K^c $ is open and therefore path-wise connected. (Proposition \ref{prop-one})\\
Let $A\supset K$. If the conclusion is true for every $A$ open and bounded and $A_1$ is an open set containing $K$, then $K\subset A_1 \cap B(0,M) \subset A_1$, for a sufficiently large $M>0$. Since $ A_1 \cap B(0,M)$ is bounded there is an open $V$ such as $K \subset V \subset A_1 \cap B(0,M) \subset A_1$ and every component of $V$ is simply connected. Therefore, we can assume, without a loss of generality, that $A$ is bounded.\\
Since $A$ is bounded, $A^c$ has exactly one unbounded component, $U_{A^c}$. (Proposition \ref{prop_two})\\
If $U_{A^c}$ is the only component of $A^c$, then $A^c$ is connected and therefore connected on the Riemann sphere (Proposition \ref{prop_three}) and consequently (Theorem \ref{theorem_one}) $A$ is simply connected and the conclusion follows.
We can, therefore, assume that $A^c$ has bounded components. \\
Since $K\subset A$ and $A$ open, $\exists \e > 0 : B(z,\e)\subset A \forall z\in K$ (i.e $\e < dist(\partial A, K)$). $K$ is compact, hence only finite such balls $B_1, ... , B_n$, $n\in\N$ with radius $\e$ suffice to cover $K$: $K \subset \bigcup\limits_{i=1}^n B_i = W \subset A$\\
Now,  $W^c$ has a finite number of components, $N\in\N$ \footnote{This holds for any open and convex subsets of $\R^n$: Let $V_1,\ldots,V_n$ be open, bounded and convex subsets of $\R^n$. Then the complement $F=\mathbb R^2\smallsetminus\bigcup_{i=1}^n V_i$ possesses only finitely many connected components.\cite{Overflow}.} If $z_1, z_2, ... , z_N \in \C$ each belong to a different component of $W^c $, then $z_1, z_2, ... , z_N \in K^c$ and since $K^c$ is path-wise connected, $K^c$ is also path-wise connected on the Riemann sphere (Proposition \ref{prop_three}). Therefore, for $ i = 1, 2, ...,N$ there exist $\g_i : [0,1] \to (\Cinf \setminus K)$  continuous, such as $\g_ i(0)=z_i$, $\g_i(0)=\infty$. Without a loss of generality, we can assume that $\g_i(s) \in \C$  $\forall s \in [0,1)$ and $\forall i \in \{1, 2, ..., N\}$\\
We now define $G_i :=\g_i([0,1))\subset \C$ for all $i \in \{1, 2, ..., N\} $ and notice that $G_i$ is closed $\forall i \in \{1, 2, ..., N\}$. Hence $V:=W\cap(\bigcup \limits _{i=1}^{N}G_i)^c\subset A$ is open and $V^c$ is path-wise connected. Therefore $V^c$ path-wise connected on Riemann sphere and hence every component of $V$ is simply connected (Corollary of \ref{cor_one}).\\

Conversely, let $K^c$ be disconnected. \\
If $U_{K^c}$ is the unbounded component of $K^c$, then there exists  a $z_0 \in \C$ such that  $z_0 \not \in K$ and  $z_0 \not \in U_{K ^c}$. Since $K^c$ is open, there exist an $\e >0 : B(z_0,\e) \subset K^c$. For this $\e$, it is straightforward to show that $B(z_,\e) \cap U_{K^c} = \varnothing $ and $\overline {B} (z_0, \frac{\e}{2}) \subset K^c$.\\
Since $K$ is compact, it is bounded and therefore there exists an  $M>0$ such that $B(z_0,M)\supset K$\\
Let $A:=\{z\in \C : \frac {\e}{2} < |z-z_0| <M\}$ the annulus centered at $z_0$, with radii $\frac {\e}{2}$ and $M$. Then there is no open subset $V$ of $A$ with $K \subset V$ such as all components of $V$ are simply connected. Indeed, for any $\g :[0,1] \to \Cinf$ continuous such as $\g(0)=z_0$ and $\g(1)=\infty$, we have $\g^* \cap K \neq \varnothing$ and since $K \subset V$, $\g^* \cap V \neq \varnothing$. Therefore $V^c$ is not path-wise connected on Riemann sphere and the conclusion follows from Theorem \ref{theorem_one}. \hfill $\Box$ 
\end{theorem} 

\noindent Our main result follows as a corollary from the characterization.

\begin{corollary}\label{cor_two}
Let $K,L \subset \C$ be disjoint and compact, such as $K^c$ and $L^c$ are connected. Then $(K\cup L)^c$ is connected.
\proof Let $A$ open such as $K\cup L \subset A$. Since $K,L$ are disjoint and closed, there are $A_1,A_2 \subset A$ open and disjoint such as $K\subset A_1$ and $L\subset A_2$. From the characterization  (theorem \ref{main} ), $\exists V_1,V_2$ open such as $K\subset V_1 \subset A_1$, $L\subset V_2 \subset A_2$ and every component of $V_1$ and $V_2$ is simply connected.\\
Then $K\cup L \subset V_1 \cup V_2 \cup A$ and every component of $V:=V_1\cup V_2$ is simply connected. Now, it is obvious that $V$ is open and since $A$ was arbitrary, the conclusion follows from the characterization  \hfill $\Box$
\end{corollary}
\newpage

%\begin{prop}\label{prop A2}
%Let $A \subset \C$ open and $\g: [0,1] \to A$. Then there exists a finite piecewise linear curve which is $A$-null-homotopic to $\g$.
%\proof Since $\g^*=\g([0,1])$ is compact, it can be covered by a finite number of open disks. In particular, there exist  $n\in\N$ and $D_1, ..., D_n$ open disks, such as $D_i \subset A$ $\forall i=1, ... , n$ and $\g^* \subset \bigcup\limits _{i=1} ^n {D_i}$.\\ The preimages $\g^{-1}(D_i)$ of these disks form an open cover of $[0,1]$ which, by the Lebesgue's number lemma, has a Lebesgue number $\rho > 0$ such that, for any partion $0=t_0<t_1< ... < t_k = 1$ of $[0,1]$ with $\max\limits_{l=1,...,k} |t_l-t_{l-1}| < \rho$, we have that for all $l=1,...,k$ there exists an $i(l)\in \{ 1,...n\}$ such that $\g([t_l-t_{l-1}]) \subset D_{i(l)}$. \\
%We then join the points $\g(t_{l-1})$ and $\g(t_l)$ by the line segment $J_l:=[\g(t_{l-1}), \g(t_l) ] \subset D_{i(l)}$. Because the disks are simply connected and subsets of $A$, each $J_l$ is $A$-homotopic to $\g([t_l-t_{l-1}])$. \\ The union of these lines segments forms a polygonal loop $\s:= \bigcup \limits_{l=1}^k J_l$ which is $A$-homotopic to the original loop $\g$.
%\end{prop}

\section*{Acknowledgements}

I express my deep sense of gratitude to my teacher, Y.Smyrlis - University of Cyprus, for his help, support and guidance. Without him, this work would have never been possible.\\
I would also like to extend my sincere gratitude to V.Nestoridis - University of Athens, for initiating the paper, suggesting the subject and providing general guidance on the relevant bibliography.\\
Finally, I am very much thankful to A.Hadjigeorgiou for the fruitful discussions on this topic. \\
  
\bibliographystyle{amsplain}

\begin{thebibliography}{2}

\bibitem{CP}  C. Chui and M.N. Parnes, Approximation by overconvergence of power series , J. Math. Anal. Appl. 36 (1971), 693–696

\bibitem{Luh}W. Luh, Approximation analytischer Funktionen durch uberkonvergente Potenzreihen und
deren Matrix- Transformierten, Mitt. Math.Sem. Giessen 88 (1970), 1–56.


\bibitem{N} V. Nestoridis, Universal Taylor series, Ann. Inst. Fourier 46 (1996), 1293–1306.

\bibitem{Forster} O.Forster, Lectures on Riemannian Surfaces, Graduate Texts in Mathematics 81, Springer 20-33

\bibitem{Rudin} W.Rudin, Real and complex analysis, 3rd ed.
McGraw-Hill, Inc. New York, NY, USA ©1987 
ISBN:0070542341, 268-269
%
%\bibitem{Munkres} James R. Munkres, Topology, Prentice Hall, Incorporated, 2000
%ISBN	0131816292, 9780131816299

%@MISC {Overflw,
%    TITLE = {Complement of a finite union of convex sets},
%    AUTHOR = {Eric Wofsey (http://mathoverflow.net/users/75/eric-wofsey)},
%    HOWPUBLISHED = {MathOverflow},
%    NOTE = {URL:http://mathoverflow.net/q/224294 (version: 2015-11-23)},
%    EPRINT = {http://mathoverflow.net/q/224294},
%    URL = {http://mathoverflow.net/q/224294}
%}
\bibitem{Overflow} Eric Wofsey (http://mathoverflow.net/users/75/eric-wofsey), Complement of a finite union of convex sets, (http://mathoverflow.net/q/224294), 2015-11-23 

\end{thebibliography}

\end{document}